\documentclass{article}
\pdfoutput=1
\usepackage[english]{babel}
\usepackage{graphicx}
\usepackage{amsmath}
\usepackage{amssymb, wasysym, amsfonts, latexsym,amsthm}
    
\newcommand{\X}{\mathfrak X}

\newcommand{\romb}{\diamondsuit}
\newcommand{\Log}{\mathop{Log}}
\newcommand{\logic}[1]{\mathsf{#1}}

\newcommand{\ML}{\mathcal{ML}}

\newcommand{\lK}{\logic{K}}

\newcommand{\Lo}{\mathsf{L}}
\newcommand{\rest}[1]{\bigl. #1 \bigr|}

\newcommand{\pmor}{\twoheadrightarrow}

\newcommand{\cM}{\mathcal{M}}
\newcommand{\cF}{\mathcal{F}}

\newcommand{\abs}[1]{\left\vert#1\right\vert}
\newcommand{\set}[1]{\left\{#1\right\}}
\newcommand{\Set}[1]{\bigl\{#1\bigr\}}
\newcommand{\setdef}[2][x]{\set{#1\,\left |\,#2 \right .}}
\newcommand{\Setdef}[2][x]{\Set{\bigl.#1\,\bigr|\,#2 }}

\newcommand{\Real}{\mathbb R}
\newcommand{\QQ}{\mathbb Q}
\newcommand{\Natr}{\mathbb N}

\newcommand{\Y}{\mathcal Y}
\newcommand{\ve}[1]{\mathbf{#1}}

\newcommand{\Nf}{\mathcal{N_\omega}}
\newcommand{\N}{\mathcal{N}}

\newcommand{\wprd}[1]{\langle #1 \rangle}

\newcommand{\len}{\mathop{len}}

\newtheorem{theorem}{Theorem}[section]

\newtheorem{proposition}[theorem]{Proposition}

\newtheorem{corollary}[theorem]{Corollary}

\newtheorem{lemma}[theorem]{Lemma}

\theoremstyle{definition}

\newtheorem{definition}[theorem]{Definition}

\theoremstyle{remark}

\newtheorem{remark}[theorem]{Remark}
\newtheorem*{remark*}{Remark}

\newcommand{\st}{\mathop{st}}

\begin{document}

  \title{On neighborhood product of some Horn axiomatizable logics.}

  \author{Andrey Kudinov\footnote{This work was partialy supported by RFBR grants N 14-01-31442-mol-a.}\\ \small
kudinov\ [at here]\ iitp\ [dot]\ ru \\
\small
Institute for Information Transmission Problems, Russian Academy of
Sciences\\ \small
National Research University Higher School of Economics, Moscow, Russia}
\date{}

 \maketitle
 
\begin{abstract}
  We consider modal logics of products of neighborhood frames. We define n-product of modal logics as the logic of all products of neighborhood frames of corresponding logics and find n-product of any two pretransitive Horn axiomatizable logics. As a corrolary we find the d-logic of products of topological spaces for some classes of topological spaces.
\end{abstract}
\smallskip 


\section{Introduction}

Neighborhood semantics is a generalization of Kripke semantics and topological semantics. It was introduced independently by Dana Scott \cite{Scott1970-SCOAOM-2} and Richard Montague \cite{montague1970}.
In this paper we consider the product of neighborhood frames introduced by Sano in \cite{Sano11:AxiomHybriProduMonotNeighFrame}. It is a generalization of the product of topological spaces\footnote{``Product of topological spaces'' is a well-known notion in Topology but it is different from what we use here (for details see \cite{benthem:MultimodalLogicsProductsTopologies})} presented in \cite{benthem:MultimodalLogicsProductsTopologies}.

The product of neighborhood frames is defined in the same manner as the product of Kripke frames (see \cite{Sege73:Twodimodallogic} and \cite{Sheh78:Twodimodallogic}). But there are some differences.
Axioms of commutativity and Church-Rosser property are valid in any product of Kripke frames. Whereas in \cite{benthem:MultimodalLogicsProductsTopologies} it was shown that the logic of the products of all topological spaces is the fusion of logics $\logic{S4} \ast \logic{S4}$. Even more, $\logic{S4} \ast \logic{S4}$ is complete w.r.t.{} the product $\mathbb{Q} \times_t \mathbb{Q}$ ($\times_t$ stands for product of topological spaces, defined in \cite{benthem:MultimodalLogicsProductsTopologies}).

In \cite{kudinov_aiml12} it was proved that for any pair $\logic{L}$ and $\logic{L'}$  of logics from $\set{\logic{S4}, \logic{D4}, \logic{D}, \logic{T}}$
modal logic of the family of products of $\logic{L}$-neighborhood frames and $\logic{L'}$-neighborhood
frames is the fusion of $\logic{L}$ and $\logic{L'}$. But at that point it was unclear how to proceed in the case of logics that do not contain axiom $\romb \top$ (correspond to seriality).
In paper \cite{kudinov2014neighbourhood} we showed that any product of neighborhood frames in fact satisfy axiom $B \to \Box_2 B$, where $B$ is a variable-free and $\Box_2$-free formula (and the similar is true for $\Box_1$). We proved that $\lK \ast \lK$ plus all such axioms is the logic of all products of neighborhood frames. For any two modal logics $ \Lo_1 $ and $ \Lo_2 $ we can define $ \wprd{\Lo_1, \Lo_2} $ as $ \Lo_1 \ast \Lo_2 $ plus all the axioms from above.

In this paper we find a sufficient conditions for two logics $ \Lo_1 $ and $ \Lo_2 $ to be n-product matching. That is that $ \Lo_1 \times_n \Lo_2 = \wprd{\Lo_1, \Lo_2} $, where $ \Lo_1 \times_n \Lo_2 $ is the logic of all products of neighborhood frames $ \X_1 \times \X_2 $ such that $ \X_1 \models \Lo_1 $ and $ \X_2 \models \Lo_2 $.

Neighborhood frames are often considered in the context of non-normal modal logics, since, unlike Kripke semantics, it is complete w.r.t.{} many non-normal logics. 
As for the normal modal logics, neighborhood frames rarely give anything new in comparison to Kripke frames. 
This paper, however, shows that in case of products normal neighborhood frames, that correspond to normal modal logics, give different results from Kripke frames. 

The results of this paper (and others: \cite{benthem:MultimodalLogicsProductsTopologies}, \cite{kudinov_aiml12}, \cite{Sano11:AxiomHybriProduMonotNeighFrame}) show that ``neighborhood'' product, in general, generate weaker logic in comparison to ``Kripke'' product. It also shows how the notion of the product of modal logics depends on the underlining semantics.

The results of this paper have corollaries for derivational semantics of topological spaces. In particular the logic of all products of all $ T_1 $ spaces is $ \wprd{\logic{K4}, \logic{K4}} $.

\section{Language, logics and semantics}
In this paper we study propositional modal logics. A formula is defined recursively as follows:
$$
A ::= p\; |\;\bot \; | \; A \to A \; | \; \Box_i A,
$$
where $p\in \mathrm{PROP}$ is a propositional letter, and $\Box_i$ is a modal operator $ i= 1, \ldots n$.
Other connectives are introduced as abbreviations: classical connectives are expressed through $\bot$ and $\to$,
dual modal operators $\Diamond_i$ are expressed as $\lnot\Box_i\lnot$. The set of all modal formulas is denoted as $ \mathcal{ML}_n $, and in order to specify the modalities used in the language we right them in subindex, for example: $ \mathcal{ML}_{\Box_1} $ or $ \mathcal{ML}_{\Box_2} $. 

\begin{definition} A \emph{normal modal logic\/} (or \emph{a logic,\/} for short) is
a set of modal formulas closed under Substitution $\left
(\frac{A(p_i)}{A(B)}\right )$, Modus Ponens $\left (\frac{A,\,
A\to B}{B}\right )$ and Generalization rules $\left
(\frac{A}{\Box_i A}\right)$, containing all
classic tautologies and normality axioms:
 $$
 \begin {array}{l}
 \Box_i (p\to q)\to (\Box_i p\to \Box_i q).
 \end{array}
$$

$\logic{K_n}$ denotes \emph{the minimal normal modal logic with $n$ modalities} and $\lK = \logic{K_1}$.
\end{definition}

Let $\logic{L}$ be a logic and let $\Gamma$ be a set of formulas, then $\logic{L} +
\Gamma$ denotes the minimal logic containing $\logic{L}$ and $\Gamma$. If
$\Gamma = \set{A}$, then we write $\logic{L} + A$ rather than $\logic{L}+  \{ A \}$.

\begin{definition}
 Formula $B$ is called \emph{closed} if it does not contain variables.
\end{definition}

\begin{definition}
Let $\logic{L_1}$ and $\logic{L_2}$ be two modal logics with one modality $\Box$, then \emph{fusion} of these logics is the following modal logic with 2 modalities:
\[
\logic{L_1} \ast \logic{L_2} = \logic{K_2} + \logic{L'_1} + \logic{L'_2};
\] 
where $\logic{L'_i}$ is the set of all formulas from $\logic{L_i}$ where all $\Box$ are replaced by $\Box_i$.
\end{definition}

\begin{definition}
 Let $R \subseteq W\times W$ be a relation on $W \ne \varnothing$, then for $k\ge 1$ and $w \in W$ we define
 \begin{align*}
  R^0 &= Id_W  = \setdef[(w,w)]{w\in W};\\
  R^{n+1} &= R^n \circ R;\\
  R^* &= \bigcup_{k=0}^{\infty} R^k;\\
  R(w) &= \setdef[u]{ w R u}.
 \end{align*}
 
A \emph{Kripke frame} with $n$ relations is a tuple $(W, R_1, \ldots R_n)$, where $ W $ is a non-empty set and $ R_i \subseteq W \times W $ is a relation on $ W $ for each $ i\in \set{1, \ldots n}$.

\begin{remark*}
	We will sometimes write $w\in F$ as a shortcut for $w \in W$ and $F = (W, R_1, \ldots R_n)$.
\end{remark*}

A frame $ F $ with a \emph{valuation} $ V:PROP \to 2^W $ is called \emph{a model} $M = (F, V)$.  

For a Kripke frame $F = (W, R_1, \ldots R_n)$ we define the \emph{subframe} generated by $w \in W$ as the frame $F^w = (W', R|_{W'})$, where $W' = (R_1 \cup \ldots \cup R_n)^*(w)$ and ${R_i}|_{W'} = R_i \cap W' \times W'$. A frame $F$ is called \emph{rooted} if $F=F^w$ for some $w$.
 
 The truth of a formula in a model $ M $ at a point $ x\in W $ is defined as usual by induction on the length of the formula:
 \begin{align*}
  M, x &\not\models \bot &\\
 M, x &\models p &\iff& x \in V(p)\\
 M, x &\models A \to B &\iff& M, x \not\models A \hbox{ or } M, x \models B\\
 M, x &\models \Box_i A &\iff& \forall y\; (xR_i y \Rightarrow M, y \models A)
 \end{align*} 
\end{definition}

A formula is \emph{valid in a Kripke model} $M$ if it is true at all points of $M$ (notation $M \models A$).
A formula is \emph{valid in an Kripke frame} $F$ if it is valid in all models based on $F$ (notation $F \models A$). 
We write $F \models \Lo$ if for any $A \in \Lo$, $F\models A$.
\emph{Logic} of a class of Kripke frames $\mathcal{C}$ is $\Log(\mathcal{C}) = \setdef[A]{F \models A \hbox{ for all } F \in \mathcal{C}}$.
For logic $\Lo$ we also define $V(\Lo) =\setdef[F]{\hbox{$F$ is an Kripke frame and } F\models \Lo}$. Note, that if there are no $F$ such that $F \models \Lo$, then $V(\Lo) = \varnothing$.

\begin{definition} \label{def:p_morphism}
	Let $F = (W, R_1, \ldots R_n)$ and $G = (U, S_1, \ldots S_n)$ be Kripke frames. Then function $f: W \to U$ is a \emph{p-morphism} if 
	\begin{enumerate}
		\item $f$ is surjective;
		\item {}\textbf{[monotonisity]} for any $w, v\in W$ from $ wR_i v $ follows $ f(w) S_i f(v) $;
		\item {}\textbf{[lifting]} for any $w\in W$ and $v' \in U$ such that $ f(w)S_i v' $ there exists $v \in W$ such that $ wRv $ and $ f(v)=v' $.
	\end{enumerate}
	In notation $f: F \pmor G$.
\end{definition}

The following is known as the p-morphism lemma

\begin{lemma}\label{lem:pmorhism4K-frame}
	Let $f: F \pmor G$ and $V$ be a valuation on $G$. We define a valuation on $F$: $\left[f^{-1}(V)\right](p) = f^{-1}(V(p))$. Then for any $ w\in F $ and formula $ A $
	\[
	F, f^{-1}(V), w \models A \iff G, V, f(w) \models A.
	\] 
\end{lemma}

The proof is by standard induction on the length of formula $A$. The following is a straightforward corollary. 

\begin{corollary}\label{cor:pmorphKframes}
	If $f: F \pmor G$, then $\Log(F) \subseteq \Log(G)$.
\end{corollary}

For a modal logic $\Lo$ with $n$ modalities we define the canonical model (cf. \cite{blackburn_modal_2002}) $\cM_{\Lo} = (\cF_{\Lo}, V_{\Lo})$, where $\cF = (W, R_1, \ldots R_n)$ such that
 \begin{align*}
  W &= \setdef[x]{x \hbox{--- is an $\Lo$-complete set of formulas}},\\
  x R_i y &\iff \forall A ( \Box A \in x \Rightarrow A \in y),\\
  x \in V(p) &\iff p \in x.
 \end{align*}

The classical result on canonical models is
\begin{lemma}
For any formula $ A $ and any logic $ \Lo $
	$$ \cM_{\Lo}, x \models A \iff A \in x.$$
\end{lemma}

We also define 0-canonical frame $ \cF^0_{\Lo} $ being the counterparts of canonical frame in the modal language without variables. More precisely
 \begin{align*}
 \cF^0 &= (W^0, R'_1, \ldots R'_n)\\
 W^0 &= \setdef[\bar x]{\bar x \hbox{  is an $\Lo$-complete set of \emph{closed} formulas}},\\
 \bar x R'_i \bar y &\iff \forall A ( \Box A \in \bar x \Rightarrow A \in \bar y),\\
 \end{align*}

Note that there are no 0-canonical models since there are no variables in closed formulas. So the lemma for canonical model transforms into

\begin{lemma}
For any closed formula $ A $ and any logic $ \Lo $
$$ \cF^0_{\Lo}, \bar x \models A \iff A \in \bar x.$$	
\end{lemma}





Now we are going to describe a construction of continuum unravelling. It is similar to the construction in \cite[Lemma 4.9]{GSh_Product_1998}.

\begin{definition}
	Let $ F=(W, R) = F^{w_0}$ be a rooted Kripke frame, $ S $ be a non-empty set and $ x_0 \in S $ be a fixed point in it. Then
	\begin{align*}
	F\cdot S &= (W\times S, R \cdot S),\\
	(w,x) R \cdot S (v,y) &\iff wRv,\\
	F_S = (F\cdot S)^{(w_0, x_0)} &= (W_S, R_S) \hbox{ --- a rooted subframe}.
	\end{align*}
	$ F_S $ is called the \emph{thickening} of $ F $ by $ S $.
\end{definition}

The proof of the following lemma is straightforward.
\begin{lemma}
	The first projection $ p_1(w,x) = w$ is a p-morphism $ p_1: F_V \pmor F $.
\end{lemma}

The following construction is well-known (c.f.{} \cite{blackburn_modal_2002}).

\begin{definition}
	Let $ F=(W, R_1, \ldots R_n) = F^{w_0}$ be a rooted Kripke frame. 
	
	We define the unravelling of it and a map $ \pi $ as follows
	\begin{align*}
	F^{\sharp} &= (W^{\sharp}, R^{\sharp}_1 \ldots R^{\sharp}_n)\\
	W^{\sharp} &= \setdef[w_0 R_{j_1} w_1 \ldots R_{j_m} w_m]{\forall i \in \set{1, \ldots m}  \left (w_{i-1} R_{j_i} w_{i}\right )}\\
	\pi (w_0 w_1 \ldots w_m) &= w_m,\\
	\alpha R^{\sharp}_j \beta &\iff \beta = \alpha w_{m+1} \hbox{ and } \pi(\alpha) R_j \pi(\beta). 
	\end{align*}
\end{definition} 

\begin{lemma}\label{lem:unrev_pmorph}
	The map $ \pi $ is a p-morphism: $ \pi: F^{\sharp} \pmor F $.
\end{lemma}

The proof is straightforward.

\begin{definition}
	Let $ F=F^{w_0} $ we define be a rooted frame, then we define the \emph{continuum unravelling} of it as $F^{\sharp}_{\Real} = (F_{\Real})^{\sharp} $ (the unravelling of the thickening by $ \Real $ with $0$ as the the fixed point).
\end{definition}


Furthermore, we consider neighborhood frames (c.f.{} \cite{Segerberg1971} and \cite{Chellas1980}).

\begin{definition}
	Let $X$ be a nonempty set, then $\mathcal{F} \subseteq 2^X$ is a \emph{filter} on $X$ if
	\begin{enumerate}
		\item $X \in \mathcal{F}$;
		\item if $U_1,\; U_2 \in \mathcal{F}$, then $U_1 \cap U_2 \in \mathcal{F}$;
		\item if $U_1 \in F$ and $U_1 \subseteq U_2$, then $U_2 \in \mathcal{F}$.
	\end{enumerate}
	Note, it is usually demanded that $\varnothing \notin \mathcal{F}$ ($\mathcal{F}$ is a proper filter), but we will not demand this in this paper.
\end{definition}

\begin{definition}
	A \emph{(normal) neighborhood frame} (or an n-frame) is a pair $\X = (X, \tau)$, where $X$ is a nonempty set and $\tau: X \to 2^{2^X}$ such that $\tau(x)$ is a filter on $X$ for any $x$. Function $\tau$ is called the \emph{neighborhood function} of $\X$, and sets from $\tau(x)$ are called \emph{neighborhoods of $x$}.
	\emph{The neighborhood model} (n-model) is a pair $(\X, \theta)$, where $\X = (X, \tau)$ is an n-frame and $\theta: PROP \to 2^X$ is a \emph{valuation}.  
	In a similar way, we define \emph{neighborhood 2-frame} (n-2-frame) as $(X, \tau_1, \tau_2)$ such that $\tau_i (x)$ is a filter on $X$ for any $x$, and a \emph{n-2-model}.
\end{definition}

\begin{remark*}
Note that usually neighorhood frames are defined without demanding anything from the neighborhood function. And many papers consider the \emph{monotone} neighorhood frames are considered, that is if we demand only point 3 from definition of the filter (the set of neighorhoods is closed under supersets).
\end{remark*}

\begin{definition}
	\emph{The valuation of a formula} $A$ at a point of an n-model $M = (\X, \theta)$ is defined by induction. For Boolean connectives the definition is usual, so we omit it. For modalities the definition is as follows:
	\[
	M, x \models \Box_i A \iff \exists U \forall y ( y \in U\in \tau_i(x) \Rightarrow M,y \models A).
	\]
	Formula is valid in an n-model $M$ if it is valid at all points of $M$ (notation $M \models A$).
	Formula is valid in an n-frame $\X$ if it is valid in all models based on $\X$ (notation $\X \models A$). 
	We write $\X \models L$ if for any $A \in L$, $\X\models A$.
	Logic of a class of n-frames $\mathcal{C}$ as $\Log(\mathcal{C}) = \setdef[A]{\X \models A \hbox{ for all } \X \in \mathcal{C}}$.
	For logic $L$ we also define $nV(L) =\setdef[\X]{\hbox{$\X$ is an n-frame and } \X\models L}$. Note, that if there is no $\X$ such that $\X \models L$, then $nV(L) = \varnothing$.
\end{definition}

\begin{definition}
	Let $F = (W, R)$ be a Kripke frame. We define n-frame $\N(F) = (W, \tau)$ in the following way
	\[
	\tau (w) = \setdef[U]{ R(w) \subseteq U \subseteq W}.
	\] 
\end{definition}

\begin{lemma}\label{lem:n-frame_from_kframe}
	Let $F = (W, R)$ be a Kripke frame. Then
	\[
	\Log(\N(F)) = \Log(F).
	\]
\end{lemma}

The proof is straightforward (see \cite{Chellas1980}). 

\begin{definition} \label{def:bounded_morphism}
	Let $\X = (X, \tau_1, \ldots)$ and $\Y = (Y, \sigma_1, \ldots)$ be n-frames. Then function $f: X \to Y$ is a \emph{p-morphism} if 
	\begin{enumerate}
		\item $f$ is surjective;
		\item for any $x\in X$ and $U \in \tau_i(x)$ $f(U) \in \sigma_i(f(x))$;
		\item for any $x\in X$ and $V \in \sigma_i(f(x))$ there exists $U \in \tau_i(x)$ such that \mbox{$f(U) \subseteq V$}.
	\end{enumerate}
	In notation $f: \X \pmor \Y$.
\end{definition}

\begin{remark}
	According to Lemma \ref{lem:n-frame_from_kframe}, a Kripke frame is a particular case of a neighborhood frame. It is easy to check that for any two Kripke frames $F$ and $G$ function $f$ is a p-morphism (Def \ref{def:p_morphism}) from $F$ to $G$ iff $f$ is a p-morphism (Def. \ref{def:bounded_morphism}) from $\N(F)$ to $\N(G)$. So, p-morphism for n-frames is a natural generalization of p-morphism for Kripke frames. This is why we use the same name for these formally different notions.
\end{remark}

\begin{lemma}\label{lem:pmorhism4n-frame}
	Let $\X = (X, \tau_1, \ldots)$, $\Y = (Y, \sigma_1, \ldots)$ be n-frames and $f: \X \pmor \Y$. Let $\theta$ be a valuation on $\Y$. We define $\left[f^{-1}(\theta)\right](p) = f^{-1}(\theta(p))$. Then
	\[
	\X, f^{-1}(\theta), x \models A \iff \Y, \theta, f(x) \models A.
	\] 
\end{lemma}

The proof is by standard induction on the length of formula $A$. The following is a straightforward corollary. 

\begin{corollary}\label{cor:pmorph}
	If $f: \X \pmor \Y$, then $\Log(\X) \subseteq \Log(\Y)$.
\end{corollary}

\section{Products: from Kripke to neighborhood frames}

\begin{definition}
 Let $F_i = (W_i, R_i)$ ($i=1,2$) be two Kripke frames. We define their \emph{product} (see \cite{GSh_Product_1998}) as a bimodal frame $F_1 \times F_2 = (W_1 \times W_2, R_1^h, R_2^v)$, where
 \begin{align*}
  (x,y) R_1^h (z,t) &\iff x R_1 z \ \&\ y=t,\\
  (x,y) R_2^v (z,t) &\iff x=z \ \&\ y R_2 t.\\
 \end{align*}
\end{definition}

\begin{definition}
 Let $\X_1=(X_1, \tau_1)$ and $\X_2 =(X_2, \tau_2)$ be two n-frames. Then \emph{the product} of these n-frames is an n-2-frame defined as follows:
\begin{align*}
\X_1 \times \X_2 &= (X_1 \times X_2, \tau_1^h, \tau_2^v),\\
\tau_1^h(x_1, x_2) &= \setdef[U\subseteq X_1 \times X_2]{\exists V( V \in \tau_1(x_1) \;\&\; V \times \set{x_2}  \subseteq U)},\\
\tau_2^v(x_1, x_2) &= \setdef[U\subseteq X_1 \times X_2]{\exists V( V \in \tau_2(x_2) \;\&\; \set{ x_1} \times V \subseteq U)}.
\end{align*}

\end{definition}

\begin{remark*}
Note that the product of n-frames is closed under superset. So it is possible to define product of any two monotone n-frames and it was done in \cite{Sano11:AxiomHybriProduMonotNeighFrame}. However, we consider only normal n-frames in this paper, because we use Kripke semantics in the proof of completeness.   
\end{remark*}

\begin{definition}
For two unimodal logics $\logic{L_1}$ and $\logic{L_2}$, so that $nV(\logic{L_1}) \ne \varnothing$ and $nV(\logic{L_2}) \ne \varnothing$, we define \emph{n-product} of them as follows:
\[
\logic{L_1} \times_n \logic{L_2} = \Log \bigl(\setdef[\X_1 \times \X_2]{\X_1\in nV(\logic{L_1}) \;\&\; \X_2 \in nV(\logic{L_2})}\bigr).
\]
\end{definition}

If we forget about one of its neighborhood functions, say $\tau'_2$, then $\X_1 \times \X_2$ will be a disjoint union of $\logic{L_1}$ n-frames. Hence,

\begin{proposition}[\cite{Sano11:AxiomHybriProduMonotNeighFrame}]\label{prop:fusin_of_n-frames}
 For two unimodal normal logics $\logic{L_1}$ and $\logic{L_2}$
\[
\logic{L_1} \ast \logic{L_2} \subseteq \logic{L_1} \times_n \logic{L_2}.
\]
\end{proposition}

In Chapter \ref{sec:seriality} we will show that n-product of any two logics from set $\set{\logic{S4}, \logic{D4}, \logic{D}, \logic{T}}$ equals to the fusion of corresponding logics. But this is not the case for $\lK$:

\begin{proposition}
 $\lK \times_n \lK \ne \lK \ast \lK$.
\end{proposition}

\begin{proof}
Let $\X_1=(X_1, \tau_1)$ and $\X_2 =(X_2, \tau_2)$ be two n-frames and $\X_1 \times \X_2 = (X_1 \times X_2, \tau_1^h, \tau_2^v)$.
 Consider formula $\Box_1 \bot \to \Box_2 \Box_1 \bot$. Since this formula has no variables, the truth of this formula does not depend on the valuation. So
 \begin{align*}
 \X_1 \times \X_2, (x,y) \models \Box_1 \bot &\iff \varnothing \in \tau_1^h(x,y)\iff \\ \varnothing \in \tau_1 (x) &\iff  \forall y'\in X_2\ (\varnothing \in \tau'_1(x,y')) \iff \\
 \forall y'\in X_2\ (\X_1 \times \X_2, (x,y') \models \Box_1 \bot)  &\ \Longrightarrow
 \X_1 \times \X_2, (x,y) \models \Box_2 \Box_1 \bot.
 \end{align*}
 
 Hence, $\X_1 \times \X_2 \models \Box_1 \bot \to \Box_2 \Box_1 \bot$.
 \end{proof}

 Moreover,
 
\begin{lemma}\label{lem:closed_formulas}
For any two n-frames $\X_1$ and $\X_2$ 
 1) if $B$ is a closed formula without $\Box_2$, then for any two n-frames $\X_1$ and $\X_2$ 
\[
\X_1 \times \X_2 \models B \to \Box_2 B,
\]
2) if $B$ is a closed formula without $\Box_1$, then 
\[
\X_1 \times \X_2 \models B \to \Box_1 B.
\]
\end{lemma}

\begin{proof}
We prove only 1) because 2) can be proved analogously.
Since $B$ does not contain neither $\Box_2$, nor variables, its value does not depend on the second coordinate. Let $F = \X_1 \times \X_2$.
So if $F, (x,y) \models B$, then $\forall y' (F, (x,y') \models B)$, hence, $F, (x,y) \models \Box_2 B$.
\end{proof}
We put
\[
\Delta = \setdef[B_1 \to \Box_2 B_1]{B_1 \mbox{ is closed and $\Box_2$-free}} \cup 
\setdef[B_2 \to \Box_1 B_2]{B_2 \mbox{ is closed and $\Box_1$-free}}.
\]
\begin{definition}
 For two unimodal logics $L_1$ and $L_2$, we define 
 $$\langle L_1, L_2 \rangle = L_1\ast L_2 + \Delta.$$ 
\end{definition}

From Lemma \ref{lem:closed_formulas} and Proposition \ref{prop:fusin_of_n-frames} follows

\begin{lemma}\label{lem:wprdnessesity}
For any two normal modal logics $L_1$ and $L_2$
 $\wprd{ L_1, L_2 } \subseteq L_1 \times_n L_2$.
\end{lemma}

\begin{corollary}\label{cor:KtimesK}
 $\langle \lK, \lK \rangle \subseteq \lK \times_n \lK$.
\end{corollary}

The proof of the converse inclusion demands some work. 


\section{Dense neighborhood frames}
To prove completeness of a logic w.r.t.\ neighborhood frames we are still going to rely on Kripke completeness.
So we need a way to construct a neighborhood frame out of a Kripke frame in such a way, that the neighborhood frame is \emph{dense}. An n-frame is called \emph{dense} if no point in it has a minimal neighborhood. This is important because otherwise n-frames will be equivalent to Kripke frames, and any product of Kripke frames satisfies the commutativity axioms and the Church-Rosser axiom. 
In order to construct such an n-frame, we introduce 

\begin{definition}
 For a frame $F=(W, R)$ with a fixed root $a_0$ 
 we define \emph{a path with stops} as a tuple $a_0 a_1 \ldots a_n$, so that $a_i \in W$ or $a_i = 0$ and after eliminating zeros each point is related to the next one by relation $R$. To be precise, a path with stops is a tuple of the following type
 \begin{align*}
 &a_0 0^{i_1} b_1 0^{i_2} \ldots 0^{i_m} b_m,\ \hbox{where } b_j \in W,\ i_j \ge 0,\ 0^i = \underbrace{00\ldots 0}_{\hbox{\small $i$ times}},\\
 \hbox{and }  &f_0(\alpha) = a_0 R b_1 R \ldots R b_m \in W^{\sharp}.
 \end{align*}
 
 We also consider infinite paths with stops that end with infinitely many zeros. We call these sequences \emph{pseudo-infinite paths (with stops)}.
 Let $W_\omega$ be the set of all pseudo-infinite paths in $W$.
\end{definition}

The function of forgetting zeros can be extended on $ W_{\omega} $  $f_0:W_\omega \to W^\sharp$ in the following way: for a pseudo-infinite path $\alpha = a_0 a_1 \ldots a_n \ldots $ we define
\begin{align*}
\st(\alpha) &= \min \setdef[N]{\forall k > N (a_k = 0)};\\
\alpha|_k &= a_1 \ldots a_k;\\ 
f_0(\alpha) &= f_0(\alpha|_{\st(\alpha)});\\
U_k(\alpha) &= \Setdef[\beta \in W_\omega]{\alpha |_m = \beta |_m \ \&\ f_0(\alpha) R^\sharp f_0(\beta), \ m = \max(k, \st(\alpha))}.
\end{align*}

\begin{lemma}\label{lem:U_m-base}
$U_k (\alpha) \subseteq U_m (\alpha)$ whenever $k \ge m$ 
\end{lemma}
\begin{proof}
Let $\beta \in U_k (\alpha)$. Since $\alpha |_k = \beta |_k$ and $k \ge m$, then $\alpha |_m = \beta |_m$. Hence, $\beta \in U_m (\alpha)$.
\end{proof}

\begin{definition}\label{def:n-frame_from_fframe}
Due to Lemma \ref{lem:U_m-base}, sets $U_n(\alpha)$ form a filter base. So we can define
\begin{align*}
\tau (\alpha) &- \hbox{the filter with base } \setdef[U_n(\alpha)]{n \in \Natr};\\
\Nf(F) &= (W_\omega,\tau) - \hbox{is \emph{a dense n-frame based on} $F$.} 
\end{align*}
\end{definition}

Frame $\Nf(F)$ is dense unlike $\N(F)$. Indeed,
\[  
\bigcap\limits_n U_n(\alpha) = \varnothing \not\in \tau(\alpha). 
\]

\begin{lemma}\label{lem:bmorphism_wframe2nframe}
 Let $F = (W, R)$ be a Kripke frame with root $a_0$, then
\[
f_0 : \Nf(F) \pmor \N(F^\sharp).
\]
\end{lemma}

\begin{proof}
From now on in this proof we will omit the subindex in $f_0$.
Since for any $b \in W$ there is a path $a_0 a_1 \ldots a_{n-1} b$ and, hence for pseudo-infinite path $\alpha = a_0  \ldots b 0^\omega \in X$, $f(\alpha) = b$ and $f$ is surjective.

Assume, that $\alpha\in W_\omega$ and $U \in \tau(\alpha)$. We have to prove that $R^\sharp(f(\alpha)) \subseteq f(U)$. 
There exists $m$ such that $U_m(\alpha) \subseteq U$, and since $f(U_m(\alpha)) = R^\sharp(f(\alpha))$, then
\[
R^\sharp(f(\alpha)) = f(U_m(\alpha)) \subseteq f(U).
\]

Assume that $\alpha\in W_\omega$ and $V$ is a neighborhood of $f(\alpha)$, i.e.\ $R^\sharp(f(\alpha)) \subseteq V$. We have to prove that there exists $U \in \tau(\alpha)$ such that $f(U) \subseteq V$. As $U$ we take $U_m(\alpha)$ for some $m \ge st(\alpha)$, then
\[
f(U_m(\alpha))=R^\sharp(f(\alpha)) \subseteq V. \qedhere
\]
\end{proof}

\begin{corollary}\label{cor:logicNf}
 For any frame $F$ $\Log(\Nf(F)) \subseteq \Log(F)$.
\end{corollary}

\begin{proof}

 It follows from Lemmas \ref{lem:n-frame_from_kframe},  \ref{lem:bmorphism_wframe2nframe}, \ref{lem:unrev_pmorph}, and Corollary \ref{cor:pmorph} that
\[
\Log(\Nf(F)) \subseteq \Log(\N(F^\sharp)) = \Log(F^\sharp) \subseteq \Log(F). \qedhere
\]
\end{proof}

Note that it could be the case that $ \Log(\Nf(F)) \ne \Log(F) $. To see that, let consider the natural numbers with ``next'' relation. It is convenient here to look at a number as a word in a one-letter alphabet: 
\[ 
G = (\set{1}^* , S),\ 1^n S 1^m \iff m=n+1.
\] 

Obviously $ G \models \romb p \to \Box p $.

Since in $ G $ every point except the root point has only one predecessor, we can identify a point and a path from the root to this point. Therefore, points of n-frame $\Nf(G)$ can be presented as infinite sequences of 0 and 1 with only zeros at the end. 
\begin{proposition}
	$ \Nf(G) \nvDash \romb p \to \Box p $
\end{proposition}
\begin{proof}
	Consider valuation $ \theta(p) = \setdef[0^{2n}10^{\omega}]{n \in \Natr} $. In any neighborhood of point $ 0^\omega $ there are points, where $ p $ is true and there are points where $ p $ is false. Hence,
	\[ 
	\Nf(G) \models \romb p \land \romb\lnot p. \qedhere
	 \]  
\end{proof}
	
It seems that formulas that restricts branching are not preserved under $ \Nf $ operation. In section \ref{sec:Horn} we define some formulas that are preserved.

\section{Weak product of Kripke frames}

In order to prove completeness w.r.t.{} n-frames, we first establish completeness w.r.t.{} special kind of Kripke frames.  For this purpose \emph{weak product} of Kripke frames were introduced in \cite{kudinov2014neighbourhood}. Here we modify this construction a little bit. This new construction isomorphic to the old one but in some respects better.   

\begin{definition}
    Let $ \Sigma $ be a non-empty finite set (\emph{alphabet}). A finite sequence of elements from $ \Sigma $ we call \emph{words}, the empty word is $\epsilon$.     
    The set of all words we define $\Sigma^*$. We will write words without brackets or commas, e.g.{} $ a_1 a_2 \ldots a_n \in \Sigma^* $. The length of a word is the number of elements in it:
    \[ 
    \len(a_1 a_2 \ldots a_n) = n, \quad \len(\epsilon)=0.
    \]
    
    We also define concatenation of words:
    \[ 
    a_1 a_2 \ldots a_n \cdot b_1 b_2 \ldots b_m = a_1 a_2 \ldots a_n  b_1 b_2 \ldots b_m
    \] 
\end{definition}


\begin{definition}
    Let $F_1 = (W_1, R_1)$ and $F_2= (W_2, R_2)$ be two Kripke frames with roots $x_0$ and $y_0$ respectively. Let $ \Sigma=W_1 \cup W_2 $ then we define functions $ p_1, p_2: \Sigma^* \to \Sigma^* $ and $ \pi: \Sigma^*\setminus \set{\epsilon} \to \Sigma$ by induction
    \begin{align*}
        p_1 (\epsilon) = \epsilon& \\
        p_2 (\epsilon) = \epsilon& \\
        p_1 (\ve a u) = p_1(\ve a)\cdot u&\hbox{ for } \ve a \in \Sigma^*,\ u \in W_1, \\
        p_1 (\ve a u) = p_1(\ve a)&\hbox{ for } \ve a \in \Sigma^*,\ u \in W_2,\\
        p_2 (\ve a u) = p_2(\ve a)&\hbox{ for } \ve a \in \Sigma^*,\ u \in W_1, \\
        p_2 (\ve a u) = p_2(\ve a)\cdot u&\hbox{ for } \ve a \in \Sigma^*,\ u \in W_2,\\
        \pi(\ve a u) = u &\hbox{ for } \ve a \in \Sigma^*,\ u \in \Sigma.\\
    \end{align*}
\end{definition}

Since $ F_1 $ and $ F_2 $ are frames with roots and have only one relation we will assume that paths in them do not contain relations and start from roots:
\begin{align*}
    W_1^\sharp &= \setdef[x_1\ldots x_n]{x_0 R_1 x_1 R_1\ldots R_1 x_n \hbox{ is a path in the usual sence}}\\    
    W_2^\sharp &= \setdef[y_1 \ldots y_n]{y_0 R_2 y_1 R_2 \ldots R_2 y_n \hbox{ is a path in the usual sence}}\\
\end{align*}

We define the \emph{entanglement} of $ F_1 $ and $ F_2 $ as follows
\begin{align*}
    F_1 \taurus  F_2 &= \setdef[\ve{a} \in \Sigma^*]{ p_1(\ve{a}) \in W_1^{\sharp} \hbox{ and } p_2(\ve{a}) \in W_2^{\sharp}}
\end{align*}

We define the \emph{weak product of frames} $ F_1 $ and $ F_2 $ as follows:
\begin{align*}
    \wprd{F_1, F_2} &= (F_1 \taurus  F_2, R_1^{<}, R_2^{<})\\
    \ve{a} R_1^{<} \ve{b} &\iff \exists u \in W_1 (\ve{b} = \ve{a}u)\\
    \ve{a} R_2^{<} \ve{b} &\iff \exists v \in W_2 (\ve{b} = \ve{a}v)
\end{align*}

\begin{proposition}
    For any two rooted frames $ F_1 $ and $ F_2 $ $ \wprd{F_1,F_2} \models \Delta_i$ ($ i=1,2 $).
\end{proposition}       
\begin{proof}
    Let $ B $ be a closed $\Box_2$-free formula and $ \wprd{F_1,F_2}, \ve{a} \models B $ then for any $ v \in W_2 $ we need to show that $\wprd{F_1,F_2}, \ve{a}v \models B$.
    Indeed, frames $ \left ((R^{<}_1)^*(\ve{a}), R^{<}_1|_{(R^{<}_1)^*(\ve{a})}\right ) $ and $ \left ((R^{<}_1)^*(\ve{a}v), R^{<}_1|_{(R^{<}_1)^*(\ve{a}v)}\right ) $ are isomorphic to $ (F_1^\sharp)^{p_1(\ve{a})} $. 
    Then, since $B$ is closed and do not contain $\Box_2$,  $\wprd{F_1,F_2}, \ve{a}v \models B$.
    
    For $ \Delta_2 $ the proof is similar.
\end{proof}

The aim of this section is to prove the following theorem:
\begin{theorem}\label{thm:wprd_completness}
 Logic $\wprd{\lK,\lK}$ is complete with respect to the class of all weak products of Kripke frames.
\end{theorem}

Let $ \cF = \cF^{x_0} $ be a rooted subframe of the canonical frame of logic $ \Lo $ with two modalities.

By $ \Upsilon $ we define all closed (variable-free) modal formulas of the modal language. For a point $ x\in \cF $ we define $ \bar{x}  = x \cap \Upsilon $. 
Then let $ \cF_0^{\bar{x}_0}$ be a rooted subframe of the  0-canonical frame.

We define $ \Upsilon_{i} $ as the set of all closed formulas in the language with only $ \Box_i $ modality.

\begin{lemma}
	Let $ \cF_0 = (\bar W, \bar R_1,\bar R_2) $ be the 0-canonical frame for logic $ \Lo $, such that $ \Delta \subset \Lo $. Then 
	\begin{align*}
	\bar{x} \bar R_1 \bar y &\Rightarrow \bar x \cap \Upsilon_2 = \bar y \cap \Upsilon_2,\\  
		\bar{x} \bar R_2 \bar y &\Rightarrow \bar x \cap \Upsilon_1 = \bar y \cap \Upsilon_1.
	\end{align*}
\end{lemma}

\begin{proof}
We prove only one half, since the other half is similar.
 
 For any $ A \in \Upsilon_2 $ 
 \[ 
 A \to \Box_1 A \in \Delta \hbox{ and  } \lnot A \to \Box_1 \lnot A \in \Delta
  \]
  So
\begin{align*}
&A \in \Upsilon_2 \cap \bar x \Rightarrow  A \to \Box_1 A \in \bar x \Rightarrow A \in \bar y,\\
&A \in \Upsilon_2 \hbox{ and } A \notin \bar x \Rightarrow  \lnot A \to \Box_1 \lnot A \in \bar x \Rightarrow \lnot A \in \bar y \Rightarrow  A \notin \bar y. \qedhere
\end{align*}
\end{proof}

By a straightforward induction we get

\begin{corollary}
	Let $ \cF_0 = (\bar{W}, \bar R_1, \bar R_2) $ be the 0-canonical frame for logic $ \Lo $, such that $ \Delta \subset \Lo $. Then 
	\begin{align*}
\bar x (\bar R_1 \cup \bar R_1^{-1})^* \bar y &\Rightarrow \bar x \cap \Upsilon_2 = \bar y \cap \Upsilon_2,\\  
\bar x (\bar R_2 \cup \bar R_2^{-1})^* \bar y &\Rightarrow \bar x \cap \Upsilon_1 = \bar y \cap \Upsilon_1.
	\end{align*}
\end{corollary}

Since any closed formula is canonical the following is true.
\begin{lemma}
Let $\logic{L_1}$ and $\logic{L_2}$ be two canonical logics. Then $\wprd{\logic{L_1}, \logic{L_2}}$ is also canonical.
\end{lemma}




\begin{lemma}\label{lem:pmor_cont_unravel}
	Let $ \Lo $ be a 2-modal logic, $\Lo_i = \setdef[A]{A \in \Lo \cap \mathcal{ML}_{\Box_i}}$ $(i= 1,2)$ be the 1-modal fragments of it, $\cF_{\Lo} = (W, R_1, R_2)$ be the canonical frame of $ \Lo $ and $ a \in \cF_{\Lo} $. Let $ F_i = \bigl((\cF^0_{\Lo_i})^{a\cap \Upsilon_i}\bigr)^{\sharp}_{\Real}$  be the continuum unravelling of rooted subframe of the 0-canonical model of logic $ \Lo_i $ with root $ a\cap \Upsilon_i $ $ (i=1,2) $.
	Then for any $\ve a_0 \in F_1$ such that $\pi(\ve a_0) = (a \cap \Upsilon_1, r)$ for some $r\in \Real$ there exist a p-morphism of 1-Kripke frames $ f: F_1^{\ve a_0} \pmor (R_1^*(a), R_1|_{R_1^*(a)})$ with the following property. 
	\[
	\forall \ve b \in F_1^{\ve a_0} \forall b \in R_1^*(a) \left( f(\ve b) = b \ \Rightarrow\ \exists l\in \Real \bigl(\pi(\ve b) = (b \cap \Upsilon_1, l)\bigr)\right ).
	\]
	The same is true for $F_2$.
\end{lemma}

\begin{proof}
	We will describe the construction only for $F_1$ because for $F_2$ it is similar.
	
	To simplify formulas we assume that $ G = (R_1^*(a), R_1|_{R_1^*(a)}) = (W, R)$ and $ F_1 = (W', R') $.
	
	Since $ F_1 $ and $G$ are rooted we can define map $ f:  F_1 \to G $ recursively.
	
	\textbf{Base:} $ f(\epsilon) = x_0 $.
	
	\textbf{Step:} Assume that $ f(\ve b) = x$, $ \pi(\ve b) = (x \cap \Upsilon_1,r) $ and $\ve c \in R'(\ve b) $. We should choose the image for $ \ve c $ from $ R(x) $.
	
	For $ y,z \in R(x) $ we define a relation
	\[ y \sim z \iff y \cap \Upsilon_1 =  z \cap \Upsilon_1. \]
	It is obviously an equivalence relation.
	Let $ U = R(x) / {\sim}$ be the quotient set of $ R(x) $ by $ \sim $.
	
	Since cardinality of each equivalence class $ [y] \in U $ is no greater than cardinality of canonical frame which is no greater than continuum; then there exists a splitting of $ \Real $ indexed by elements of $ [y] $ into sets of continuum cardinality:
	\[ 
	\Real = \bigsqcup\limits_{ z \in [y]} V^{[y]}_{z} \ \hbox{and for each }\ z \in [y] \ 
	\abs{V^{[y]}_{z}} = \abs{\Real}.
	\] 
	This is due to the standard result of Set Theory: $ \abs{\Real \times \Real} = \abs{\Real} $.
	
	For a fixed $ \ve c \in R'(\ve b) $ there exists $ y\in R(x) $ and $ r' $ such that
	\[
	\pi(\ve c) = ( y \cap \Upsilon_1, r'),\ r' \in V^{[y]}_y. 
	\]
	We define
	\[ 
	f(\ve c) = y.
	\]
	
	Each point in $ (\cF^0)^{\sharp}_{\Real} $ is reachable from $ (\bar x_0, 0) $ in finitely many steps. A point reachable in $ m $ steps will appear on $ m $-th iteration. So function $ f $ is defined correctly.
	
	Let us check that $ f $ is a p-morphism.
	
	\textbf{Monotonicity.} Is obvious from construction.
	
	\textbf{Lifting.} Let $ f(\ve a) = x $ and $ x R y $ then for any $ r' \in V^{[y]}_y $ and $\ve b = \ve a R (\bar y, r')$ we have $ \ve a R' \ve b $ and $ f(\ve b) = y $.
	
	\textbf{Surjectivity}.  Since $ F_1 $ and $ G $ are rooted, and root maps to root, surjectivity follows from the lifting property. 
\end{proof}

\begin{lemma}\label{lem:wprdPmorphism}
	Let $ \Lo_1 $ and $ \Lo_2 $ be two unimodal logics and $ \cF = \cF^{a_0} $ be the rooted subframe of the canonical frame for logic $ \wprd{\Lo_1, \Lo_2} $; then there exist two rooted frames $ F_1 $, $ F_2 $ and a p-morphism $f:\wprd{F_1, F_2} \pmor \cF$. 
\end{lemma}

\begin{proof}
As $ F_1 $ and $ F_2 $ we take $ \bigl((\cF^0_{\Lo_1})^{\bar x_0}\bigr)^{\sharp}_{\Real} $ and $ \bigl((\cF^0_{\Lo_2})^{\bar y_0}\bigr)^{\sharp}_{\Real} $ respectively, where $ \bar x_0 = a_0 \cap \Upsilon_1$ and $ \bar y_0 = a_0 \cap \Upsilon_2$. Let $F_1 = (W_1, R'_1)$, $F_2 = (W_2, R'_2)$, and $\cF = (W, R_1, R_2)$.
	
	
Using Lemma \ref{lem:pmor_cont_unravel} for each $ \ve a \in \wprd{F_1, F_2} $ we fix two p-morphisms:
\begin{align*}
g^{\ve a}_1 &: F_1^{p_1(\ve a)} \pmor \left(R_1^*(a), R_1|_{R_1(a)}\right),\\
g^{\ve a}_2 &: F_2^{p_2(\ve a)} \pmor \left(R_2^*(a), R_2|_{R_2(a)}\right).
\end{align*}
	
We also make sure that they are coordinated in the following way
\[
\ve a R^{<}_i \ve b R^{<}_i \ve c \Longrightarrow g^{\ve a}_i(\ve c) = g^{\ve b}_i(\ve c),
\]
where $R^{<}_1$ and $R^{<}_2$ are the the 1st and the 2nd relations in $\wprd{F_1, F_2}$.
	
We can do it because the restriction of a p-morphism to a rooted submodel is a p-morphism.
	
Let us define a map $f:\wprd{F_1, F_2} \to \cF$ recursively. The root of $\wprd{F_1, F_2}$ maps to the root of $\cF$:
\[
f(\epsilon) = a_0.
\]
	
Assume that for $\ve a \in F_1 \taurus F_2$ the map is defined. Let $\ve b = \ve a u$. If $u\in W_i$, then
\[
f(\ve{b}) = g_i^{\ve{a}}(\ve b).
\]

Let us check that $f$ is a p-morphism. The monotonisity is due to monotonisity of $g^{\ve a}_i$. To check the lifting assume that $f(\ve a) = a$ and $a R_i b$. Then $b \in R_i(a)$ and due to surjectivity of $g^{\ve a}_i$ there exist $\ve b$ such that $\ve a R'_i \ve b $ and $f(\ve b) = b$. The surjectivity follows from the rootedness of frames, and the lifting property.	
\end{proof}

To proof of Theorem \ref{thm:wprd_completness} assume that formula $A$ is not in logic $\wprd{\lK, \lK}$ then it is refutable in a rooted subframe of the canonical frame $\cF_{\wprd{\lK, \lK}}$. By Lemma \ref{lem:wprdPmorphism} there exist $F_1$ and $F_2$  such that $\wprd{F_1,F_2}$ is a p-morphic preimage of the subframe. Hence by p-morphism lemma $A$ is refutable in $\wprd{F_1,F_2}$.

\section{N-product completeness theorem for $\wprd{\lK,\lK}$}

Let $F_1 = (W_1, R_1) = F_1^{r_1}$ and $F_2 = (W_2, R_2) = F_2^{r_2}$ be two rooted frames. Assume that $W_1 \cap W_2 = \varnothing$. Consider the product of n-frames $\X_1 = (X_1, \tau_1) = \Nf(F_1)$ and $\X_2 = (X_2, \tau_2) = \Nf(F_2)$
\[
\X = (X_1 \times X_2, \tau_1', \tau_2') = \Nf(F_1) \times \Nf(F_2).
\]

We define function $g: \X_1 \times \X_2 \to \wprd{F_1,F_2}$ by induction, as follows.

Let $(\alpha, \beta) \in \X_1 \times \X_2$, so that $\alpha = x_1 x_2 \ldots$ and $\beta = y_1 y_2 \ldots$, $x_i \in W_1 \cup \set{0}$, $y_j \in W_2 \cup \set{0}$.  
We  define $g(\alpha, \beta)$ to be the finite sequence that we get after eliminating all zeros from the infinite sequence $x_1y_1x_2y_2 \ldots$.

\begin{lemma}\label{lem:main}
    Function $g$ defined above is a p-morphism: \[g: \X_1 \times \X_2 \pmor \N\bigl(\wprd{F_1, F_2}\bigr).\]
\end{lemma}
\begin{proof} 
First we need to check that for any $ \alpha \in \Nf(F_1)$ and any $ \beta \in \Nf(F_2)$ $ g(\alpha, \beta) \in F_1 \taurus F_2$. It follows from the equalities:
\[ 
p_1(g(\alpha, \beta)) = p_1 (f_0(\alpha)), \ p_2(g(\alpha, \beta)) = p_2 (f_0(\beta)).
\]

To prove surjectivity, we take 
$\ve{z} = z_1\ldots z_n \in F_1\taurus F_2$. For $i\le n$ we define
\[
x_i = \left\{
\begin{array}{ll}
z_i, &\hbox{ if $z_i \in W_1$};\\
0, &\hbox{ if $z_i \in W_2$};  
\end{array}
\right.\qquad
y_i = \left\{
\begin{array}{ll}
0, &\hbox{ if $z_i \in W_1$};\\
z_i, &\hbox{ if $z_i \in W_2$}.  
\end{array}
\right.
\]
Let $\alpha = x_1 x_2 \ldots x_n 0^\omega$ and $\beta = y_1 y_2 \ldots y_n 0^\omega$, then $g(\alpha, \beta) = \ve{z}$.
Hence, $g$ is surjective.

The next two conditions we check only for $\tau_1$, since for $\tau_2$ it is similar.
Assume that $(\alpha, \beta) \in X_1 \times X_2$ and $U \in \tau_1(\alpha, \beta)$. We need to prove that $R_1^<(g(\alpha, \beta)) \subseteq g(U)$.
There exist $m > \max\set{st(\alpha), st(\beta)}$ such that $U_m(\alpha) \times \set{\beta} \subseteq U$ and, since $g(U_m(\alpha) \times \set{\beta}) = R_1^<(g(\alpha, \beta))$, then
\[
R_1^<(g(\alpha, \beta)) = g(U_m(\alpha)\times \set{\beta}) \subseteq g(U);
\]
where $U_m(\alpha)$ is the corresponding neighborhood from $\X_1$.

Assume that  $(\alpha, \beta) \in X_1 \times X_2$ and $R_1^<(g(\alpha, \beta)) \subseteq V$. We need to prove that there exists $U \in \tau_1(\alpha, \beta)$ such that $g(U) \subseteq V$. As $U$ we take $U_m(\alpha)\times \set{\beta}$ for some $m  > \max\set{st(\alpha), st(\beta)}$, then
\[
g(U_m (\alpha) \times \set{\beta}) = R_1^<(g(\alpha, \beta)) \subseteq V.\qedhere
\]
\end{proof}

\begin{corollary}\label{cor:logicNftimesNf}
 Let $F_1 = (W_1, R_1)$ and $F_2 = (W_2, R_2)$, then 
 $\Log(\Nf(F_1) \times \Nf(F_2)) \subseteq \Log(\wprd{F_1, F_2})$.
\end{corollary}

It immediately follows from Lemma \ref{lem:main} and Corollary \ref{cor:pmorph}.

\begin{theorem}\label{thm:main}
 Logic $\wprd{\lK, \lK}$ is complete with respect to products of normal neighborhood frames, i.e.{} 
 \begin{equation}\label{eq:main}
  \wprd{\lK, \lK} = \lK \times_n \lK.
 \end{equation}
\end{theorem}

\begin{proof}
The inclusion from left to rignt of (\ref{eq:main}) was proved in Corollary \ref{cor:KtimesK}.

The converse inclusion follows from Theorem \ref{thm:wprd_completness} and Corollary \ref{cor:logicNftimesNf}. Indeed,
\begin{align*}
\lK \times_n \lK &= \bigcap_{\X_1, \X_2 \in nV(\lK)} \Log(\X_1 \times \X_2) \subseteq \\
&\subseteq \bigcap_{F_1, F_2 - \text{Kripke frames}} \Log(\Nf(F_1) \times \Nf(F_2)) \subseteq \\
&\subseteq \bigcap_{F_1, F_2 - \text{Kripke frames}} \Log(\wprd{F_1, F_2}) \subseteq \wprd{\lK, \lK}. \qedhere
\end{align*}
\end{proof}

\section{Horn axioms}\label{sec:Horn}

\begin{definition}    
Following \cite{GSh_Product_1998}, we define \emph{universal strict Horn sentence} as a first order closed formula of the form
\[ 
\forall x \forall y\forall z_1 \ldots \forall z_n \bigl(\phi(x,y, z_1, \ldots, z_n) \to \psi(x,y)\bigr),
\]
where $ \phi(x,y, z_1, \ldots, z_n) $ is quantifier-free
positive (i.e.{} it is built from atomic formulas by using $ \land $ and $ \lor $) and $ \psi(x,y) $ is an atomic formula in signature $ \Sigma = \left\langle R^2_1, \ldots, R^2_m\right \rangle $, where $ R^2_i $ is correspond to relation $ R_i $.
\end{definition}



\begin{definition}
    A logic $ \Lo $ is called \emph{HTC-logic} (from Horn preTransitive Closed logic) if it can be axiomatized by closed formulas and formulas of the type $ \Box p \to \Box^n p $, $ n\ge 0 $. These formulas correspond to universal strict Horn sentences (see \cite{GSh_Product_1998}). 
\end{definition}

Let $ \Gamma $ be a set of universal strict Horn formulas and $F$ be a Kripke frame. By $ F^{\Gamma} $ we define the $ \Gamma $-closure of $ F $, that is the minimal (in terms of inclusion of relations) frame such that all formulas from $ \Gamma $ are valid in it. Such frame exists due to \cite{GSh_Product_1998}:

\begin{lemma}[{\cite{GSh_Product_1998} Prop 7.9}]\label{lem:HornClosure}
	For any Kripke frame $ F=(W, R_1, \ldots, R_n) $ and set of universal strict Horn formulas $\Gamma$ there exist $ F^{\Gamma} = (W, R_1^{\Gamma}, \ldots, R_n^{\Gamma}) $ such that
	\begin{itemize}
		\item $ R_i \subseteq R_i^{\Gamma} $ for all $ i\in \set{1, \ldots,n } $;
		\item $ F^{\Gamma} \models \Gamma $
		\item if $ G \models \Gamma $ and $ f:F \pmor G $ then $ f: F^{\Gamma} \pmor G$.
	\end{itemize}
\end{lemma}


\begin{definition}
 Let $\Gamma$ be a set of universal strict Horn formulas, $F = (W,R)$ be a rooted frame, $\alpha \in W_{\omega}$ and $f_0:W_{\omega} \to W^{\sharp}$ be the ``forgetting zeros'' function, then we define 
 \begin{align*}
  U_k^{\Gamma}(\alpha) &= \Setdef[\beta \in W_\omega]{\alpha |_m = \beta |_m \ \&\ f_0(\alpha) (R^\sharp)^{\Gamma} f_0(\beta), \ m = \max(k, \st(\alpha))},\\ 
\tau^{\Gamma}(\alpha) &= \setdef[V]{\exists k \left ( U_k^{\Gamma}(\alpha) \subseteq V \right )},\\
\N^{\Gamma}_{\omega}(F) &= (W_{\omega}, \tau^{\Gamma}).
 \end{align*}
\end{definition}

We also need the following obvious lemma:
\begin{lemma}
For any closed modal formula $ A $ and a p-morphism of Kripke frames $f:F \pmor G$ 
\[ 
F,x \models A \iff G, f(x) \models A.
\]
\end{lemma}

And its neighborhood analog:

\begin{lemma}
	For any closed modal formula $ A $ and a p-morphism of n-frames $f:\X \pmor \Y$ 
	\[ 
	\X,x \models A \iff \Y, f(x) \models A.
	\]
\end{lemma}  

\begin{definition}
A logic $ \Lo $ is called \emph{HTC-logic} (from Horn preTransitive \& Closed logic) if it can be axiomatized by closed formulas and formulas of the type $\Box p \to \Box^n p $. The later formula corresponds to the following condition: $ R^n \subseteq R $; and, obviously, corresponds to a universal strict Horn formula. 
\end{definition} 

In \cite{GSh_Product_1998} product matching was proved for a large class of Horn axiomatizable logics, including $ \logic{S5} $. But in our case, $\logic{S5}\times_n \logic{S5} \ne \wprd{\logic{S5}, \logic{S5}}$. In fact, since neighborhood frames correspond to topological spaces in case of transitive and reflexive logics, and due to \cite{kremertopological},
$$ \logic{S5}\times_n \logic{S5} = \logic{S5}\times \logic{S5} = [\logic{S5}, \logic{S5}] = \logic{S5} * \logic{S5} + com_{12}+com_{21}+chr.$$



\begin{lemma}
    Let $\Lo$ be a HTC logic, $\Gamma$ be the corresponding set of Horn formulas, and $F\models \Lo$. If $ \Box p \to \Box^n p \in \Lo $, then  
    \[  
    \N^{\Gamma}_{\omega}(F) \models \Box p \to \Box^n p.
    \] 
\end{lemma}

\begin{proof}  
    Let $ M = (\N^{\Gamma}_{\omega}(F), \theta)$ be a neighborhood model.
 We assume that $M, \alpha \not\models \Box^n p$, and then prove that
 $M, \alpha \not\models \Box p$, i.e.
 \[ 
 \forall m \exists \beta \in U^{\Gamma}_m(\alpha)(\beta \not\models p).
  \]
  Let us fix $m$. Then
    \begin{align*}
        &\exists  \alpha_1\in U^{\Gamma}_{m}(\alpha) \left( 
       \alpha_1 \not\models \Box^{n-1} p \right) \Rightarrow \\
\Rightarrow &\exists k_2 \alpha_2\in U^{\Gamma}_{m}(\alpha_1) \left( 
\alpha_2 \not\models \Box^{n-2} p \right) \Rightarrow\\
\qquad&\vdots\qquad\vdots\qquad\vdots\qquad\vdots\qquad\vdots\qquad\vdots\qquad\vdots \\
\Rightarrow &\exists \alpha_n\in U^{\Gamma}_m(\alpha_{n-1}) \left( 
\alpha_n \not\models p \right). \\
    \end{align*}
     By definition of $ U^{\Gamma}_m(\alpha)$
     \[ 
     f_0(\alpha) (R^\sharp)^{\Gamma} f_0(\alpha_1) (R^\sharp)^{\Gamma} \ldots (R^\sharp)^{\Gamma} f_0(\alpha_n)
     \]
     and
     \[ 
     \rest{\alpha}_m = \rest{\alpha_1}_m = \ldots  =\rest{\alpha_n}_m.
     \]
     Since $ \left (W^{\sharp}, (R^\sharp)^{\Gamma}\right ) \models \Box p \to \Box^n p$, then 
     \[ 
     f_0(\alpha) (R^\sharp)^{\Gamma} f_0(\alpha_n).
      \]
     It follows, that $ \alpha_n \in U^{\Gamma}_{m}(\alpha) $.
\end{proof}

\begin{lemma}
    Let $\Lo$ be a HTC logic, $\Gamma$ be the corresponding set of Horn formulas, and $F\models \Lo$. Then
    \[ 
    f_0: \N^{\Gamma}_{\omega}(F) \pmor \N(F^{\sharp\Gamma}).
    \]    
\end{lemma}
\begin{proof}
From now on in this proof we will omit the subindex in $f_0$.
The surjectivity was established in Lemma \ref{lem:bmorphism_wframe2nframe}.

Assume, that $\alpha\in W_\omega$ and $U \in \tau^\Gamma(\alpha)$. We need to prove that $R^{\sharp\Gamma}(f(\alpha)) \subseteq f(U)$. 
There exists $m$ such that $U^{\Gamma}_m(\alpha) \subseteq U$, and since $f(U^\Gamma_m(\alpha)) = R^{\sharp\Gamma}(f(\alpha))$, then
\[
R^{\sharp\Gamma}(f(\alpha)) = f(U^\Gamma_m(\alpha)) \subseteq f(U).
\]

Assume that $\alpha\in W_\omega$ and $V$ is a neighborhood of $f(\alpha)$, i.e.\ $R^{\sharp\Gamma}(f(\alpha)) \subseteq V$. We need to prove that there exists $U \in \tau^\Gamma(\alpha)$ such that $f(U) \subseteq V$. As $U$ we take $U^{\Gamma}_m(\alpha)$ for some $m \ge st(\alpha)$, then
\[
f(U^\Gamma_m(\alpha))=R^{\sharp\Gamma}(f(\alpha)) \subseteq V. \qedhere
\]
\end{proof}

\begin{corollary}\label{cor:HTClogicNec}
 Let $\Lo$ be a HTC-logic and $F \models \Lo$; then $\N^\Gamma_\omega(F) \models \Lo$.
\end{corollary}

\begin{lemma}\label{lem:F1F2GammaPmorphism}
Let $F_1$ and $F_2$ be two frames, $\Gamma_1$ and $\Gamma_2$ be two sets of Horn sentences corresponded to HTC-logics, then
\[
\N^{\Gamma_1}_\omega(F_1) \times \N^{\Gamma_2}_\omega(F_2) \pmor \N(\wprd{F_1, F_2}^{\Gamma_1 \cup \Gamma_2}).
\]
\end{lemma}
The proof is similar to Lemma \ref{lem:main}. The underlining sets are the same and we can take the same function $g$. So, surjectivity follows. Monotonicity and lifting are proved similar.

\begin{theorem}\label{thm:main}
Let $ \Lo_1 $\,and $ \Lo_2 $\,be two HTC-logics then 
	\[ \Lo_1 \times_n \Lo_2 = \wprd{\Lo_1, \Lo_2}. \]
\end{theorem}
\begin{proof}
By Lemma \ref{lem:wprdnessesity} $\wprd{\Lo_1, \Lo_2} \subseteq \Lo_1 \times_n \Lo_2$. 

 Let $\Gamma_1$ and $\Gamma_2$ be the sets of Horn sentences corresponding to $\Lo_1$ and $\Lo_2$. Let $A \notin \wprd{\Lo_1, \Lo_2}$; then there is a rooted subframe $\mathcal{F}$ of the canonical frame of logic $\wprd{\Lo_1, \Lo_2}$ such that $\mathcal{F} \not\models A$. Then by Lemma \ref{lem:wprdPmorphism} there are frames $F_1$ and  $F_2$ such that 
 \[
 \wprd{F_1, F_2} \pmor \mathcal{F}.
 \]
 Since $\Lo_1$, $\Lo_2$ and $\wprd{\Lo_1, \Lo_2}$ are canonical then 
 \[
 \wprd{F_1, F_2}^{\Gamma_1 \cup \Gamma_2} \pmor \mathcal{F}.
 \]
 By Lemma \ref{lem:F1F2GammaPmorphism}
\[
\N^{\Gamma_1}_\omega(F_1) \times \N^{\Gamma_2}_\omega(F_2) \pmor \N\left (\wprd{F_1, F_2}^{\Gamma_1 \cup \Gamma_2}\right ).
\] 
 By Corollary \ref{cor:HTClogicNec}
 \[
 \N^{\Gamma_1}_\omega(F_1) \models \Lo_1 \hbox{ and }\N^{\Gamma_2}_\omega(F_2) \models \Lo_2.
 \]
 At the same time 
 \[
 \N^{\Gamma_1}_\omega(F_1) \times \N^{\Gamma_2}_\omega(F_2) \not\models A.
 \]
 So $\Lo_1 \times_n \Lo_2 \subseteq \wprd{\Lo_1, \Lo_2}$. 
\end{proof}

\section{Seriality axiom}\label{sec:seriality}

Consider the seriality axiom $ \lnot \Box \bot $. By induction on the length of a formula, one can easily prove 

\begin{lemma}
If $ \lnot \Box \bot \in \Lo$ then any closed formula is $ \Lo $-equivalent to $ \bot $ or $ \top $.
\end{lemma}

The base is obvious and the step follows from 
\[ 
\vdash \Box \top \leftrightarrow \top, \qquad \lnot \Box \bot \vdash \Box \bot \leftrightarrow \bot.
 \]

\begin{lemma}\label{lem:seriality}
	For a bimodal logic $ \Lo $ if $ \Lo \vdash \lnot \Box_1 \bot$ then $ \Lo \vdash B \to \Box_2 B $ for any closed formula $ B \in \ML_{\Box_1} $.
\end{lemma}

It is a simple exercise.

\begin{corollary}
	If $ \Lo_1 $ and $ \Lo_2 $ are HTC-logics and $ \lnot \Box \bot \in \Lo_1 $,  $ \lnot \Box \bot \in \Lo_2 $ then
	\[ 
	\wprd{\Lo_1, \Lo_2} = \Lo_1 * \Lo_2.
	\]
\end{corollary}
From the above and Theorem \ref{thm:main} it follows

\begin{theorem}\label{thm:main_seriality}
Let  $\logic{L_1}$ and $\logic{L_2}$ be HTC-logics with seriality then
\[
\logic{L_1} \times_n \logic{L_2} = \logic{L_1} * \logic{L_2}.
\]
\end{theorem}

Note that this theorem covers results from \cite{kudinov_aiml12}, since logics $\logic{D}$, $\logic{T}$, $\logic{D4}$ and $\logic{S4}$ are all HTC-logics with seriality.

\begin{proposition}
    If $ \Lo_1 $ and $ \Lo_2 $ are finitely axiomatizable, and have only finitely many inequivalent closed formulas then $ \wprd{\Lo_1, \Lo_2} $ is finitely axiomatizable. 
\end{proposition}

This is because in this case $ \Delta $ has only have only finitely many inequivalent formulas. For example it is true if $ \Box^n\bot \in \Lo_1 $ and $ \Box^m\bot \in \Lo_2 $.

\section{Derivational semantics}\label{sec:derivational_semantics}

The derivational semantics studied by many authors (see for example \cite{shehtman1990derived} or \cite{bezhanishvili2005some}) can be equivalently defined as follows.

Let $ \X = (X, T) $ be a topological space. We define $ \tau^{\X}_d(x) = \setdef[U]{U' \setminus \set{x} \subseteq U,\ x\in U'\in T}$. Then for valuation $ V $ on $ X $ the following is true
\[  
\X, V, x \models_d A \iff (X, \tau^{\X}_d), V, x \models_n A.
\]
where $ \models_d $ corresponds to derivational semantics, and $ \models_n $ corresponds to neighborhood semantics.
We define $ \N_d(\X) = (X, \tau^{\X}_d) $.

For a class of topological spaces $ \mathcal{C} $ and logics $ \Lo_1 $ and $ \Lo_2 $ we put
\begin{align*}
\mathop{Log_d}(\mathcal{C}) &= \setdef[A]{\forall \X\in \mathcal C(\X \models_d A)}\\
\Lo_1 \times_{d} \Lo_2 &= \mathop{Log_d}(\setdef[\X_1 \times \X_2]{\hbox{$ \X_1 $, $ \X_2 $ -- topological spaces, }\X_1 \models_d \Lo_1, \X_2 \models_d \Lo_2})
\end{align*}

We say that $ \mathop{Log_d}(\mathcal{C})$ is the \emph{d-logic} of $ \mathcal{C} $.

\begin{theorem}
\begin{enumerate}
	\item $ \logic{K4} \times_{d} \logic{K4} = \wprd{\logic{K4}, \logic{K4}} $;
	\item $ \logic{K4} \times_{d} \logic{D4} = \wprd{\logic{K4}, \logic{D4}} $;
	\item $ \logic{D4} \times_{d} \logic{K4} = \wprd{\logic{D4}, \logic{K4}} $;
	\item $ \logic{D4} \times_{d} \logic{D4} = \logic{D4} * \logic{D4} $.
\end{enumerate}
\end{theorem}
	
\begin{proof}
 It follows from Theorems \ref{thm:main} and \ref{thm:main_seriality}.
 But, it is not a straightforward corollary, because for a logic $ \Lo $ the set of $ \Lo $-n-frames and the set of all n-frames that correspond to $ \Lo $-topological spaces do not coincide. Indeed in a topological space $ X $ family of neighborhoods of point $ x $ always contains set $ X\setminus \set{x} $ and it is not the case for n-frames.
 
 So to prove this theorem it is sufficient to say that all logics mentioned in this theorem are not reflexive and the unraveling are irreflexive. So let $ F=(W,R) $, then $ F^\sharp $ is irreflexive, and $ \N_\omega^\Gamma(F^\sharp) $ can be obtained as $ \N_d(\X) $. Where $ \X = (W_\omega, T) $, $ \Gamma $ is the Horn sentence expressing transitivity, and sets $ U^{\Gamma}_n $ form the base for topology $ T $. 
\end{proof}
    
\begin{theorem}
	The d-logic of all products of all $ T_1 $ spaces is $ \wprd{\logic{K4}, \logic{K4}} $.
\end{theorem}	

It is enough to check that topological space corresponding to $ \N_{\omega}^{\Gamma}(F) $ is a $ T_1 $ space, whenever $ F $ is the unraveling of a rooted $ \logic{S4} $-frame and $ \Gamma $ corresponds to transitivity. This can be easily checked.  

\section{Conclusions}

We are still in the beginning of the road of studying products of neighborhood frames.

This topic can be interesting from different points of view. It is interesting by itself because it is a natural way to combine modal logics, and the result is weaker then product of logics based on Kripke semantics. It is also interesting because using products we can express new properties, for example $\QQ$ and $\Real$ are indistinguishable in the unimodal language with topological semantics, whereas logics of $\QQ \times_t \QQ$ and $\Real \times_t \Real$ are different (see \cite{Kremer2015}). It is also possible that this construction will be useful for epistemic modal logic as semantics for multi-agents systems.

There are a lot of open questions in this area, to name a few:
\begin{itemize}
    \item find other sufficient conditions for product matching; 
    \item investigate products of type $ \Lo \times_n \logic{S5} $; in a forthcoming paper we will find logics of this type for any HTC logic $\Lo$ (this result is announced at AiML'16 conference);
    \item find logics of $\Real \times_t \Real$ and $\mathcal{C} \times_t \mathcal{C}$, where $\mathcal{C}$ is the Cantor space;
    \item find n-products of well known logics like $\logic{S4.1}$, $\logic{S4.2}$, $\logic{S4.3}$, $\logic{GL}$, $\logic{Grz}$, $\logic{DL}$ and other.
\end{itemize}

\bibliographystyle{plain}
\bibliography{kudinov}

\end{document}